\documentclass{conm-p-l}
\copyrightinfo{2000}{American Mathematical Society}

\newtheorem{theorem}{Theorem}[section]

\theoremstyle{definition}

\theoremstyle{remark}
\newtheorem{remark}[theorem]{Remark}

\numberwithin{equation}{section}

\newcommand{\abs}[1]{\lvert#1\rvert}
\newcommand{\Czpol}{\mathbb C_{\mathbb Z}^{\mathrm{pol}}} 
\DeclareMathOperator{\ad}{ad}
\DeclareMathOperator{\CToda}{\mathcal{C}Toda}
\DeclareMathOperator{\End}{End}
\DeclareMathOperator{\Hol}{Hol}
\newcommand{\inv}{^{-1}}
\newcommand{\KT}{K_T}
\DeclareMathOperator{\LToda}{\mathcal{L}Toda}
\newcommand{\nop}[1]{:\!#1\!:}
\DeclareMathOperator{\Sing}{Sing}
\DeclareMathOperator{\Tr}{Tr}
\DeclareMathOperator{\Res}{Res}
\DeclareMathOperator{\VToda}{VToda}

\begin{document}

\title[$H_T$-Vertex Algebras]{$H_T$-Vertex Algebras}
\author{Maarten J. Bergvelt}
\address{Department of Mathematics\\ University of Illinois\\
  Urbana-Champaign\\ Illinois 61801}
\email{bergv@uiuc.edu}

\date{\today}

\begin{abstract}
  The usual vertex algebras have as underlying symmetry the Hopf
  algebra $H_D=\mathbb C[D]$ of infinitesimal translations. We show
  that it is possible to replace $H_D$ by another symmetry algebra
  $H_T=\mathbb C[T,T\inv]$, the group algebra of the Abelian group
  generated by $T$. $H_T$ is the algebra of symmetries of a lattice of
  rank 1, and the construction gives a class of vertex algebras related to
  the Infinite Toda Lattice in the same way as the usual $H_D$-vertex
  algebras are related to Korteweg-de Vries hierarchies.
\end{abstract}

\maketitle

\section{Introduction}
\label{sec:Intro}

A vertex algebra is, roughly speaking, a \emph{singular}, commutative,
associative, unital algebra with \emph{symmetry}, see for instance
\cite{MR1865087} for a proposal to make this vague statement
precise.

The usual vertex algebras
(\cite{MR843307},\cite{MR996026},\cite{MR1651389}) have as
symmetry algebra the Hopf algebra $H_D=\mathbb C[D]$ of infinitesimal
translations ($D=\frac{d}{dz}$, say). The singularities are obtained
by localizing the dual $H_D^*=\mathbb C[[z]]$ by inverting $z$, to
obtain the Laurent series algebra $K_D=\mathbb C[[z]][z\inv]$. To
emphasize the r\^ole of $H_D$ as symmetry algebra we will refer to the
usual vertex algebras as $H_D$-vertex algebras.

We can think of $H_D$ as the universal enveloping algebra of the 1
dimensional (Abelian) Lie algebra generated by $D$. In this paper I
will sketch what a vertex algebra looks like when we take not $H_D$
but $H_T$ as symmetry algebra, where $H_T=\mathbb C[T,T\inv]$ is the
group algebra of the free Abelian group generated by $T$. More details
can be found in \cite{math.QA/0505289}.

A motivation to study $H_T$-vertex algebras comes from the theory of
integrable systems. There are many similarities between Korteweg-de
Vries type hierarchies on the one hand, and infinite Toda lattice
hierarchies on the other. For instance, both have infinitely many
conservation laws, multiple Hamiltonian structures, tau-functions and
connections with infinite dimensional Grassmann manifold, and Miura
transformations. (See \cite{MR785802} for an overview of lattice
hierarchies). It is well known that Korteweg-de Vries type hierarchies
are intimately related to $H_D$-vertex algebras, see for instance
\cite{MR2082709}. Now the infinite Toda lattice hierarchy has
$H_T$-symmetry, with the generator $T$ acting as a shift by one step
on the lattice. So it is natural to look for a type of vertex algebras
that are related to the infinite Toda lattice in the same way as
$H_D$-vertex algebras are related to Korteweg-de Vries hierarchies.

\section{Review of some ingredients for $H_D$-vertex algebras}
\label{sec:RevH_DVertexalg}

We will review the ingredients that go into the definition of
$H_D$-vertex algebras from a slightly unorthodox point of view,
in order to prepare the way for replacing $H_D$ as the symmetry
algebra by $H_T$. See also
\cite{MR1653021} or \cite{math.QA/9904104}.

Let $W_1,W_2$ be vector spaces. A \emph{$W_2$-valued distribution} on $W_1$
is just a linear map $\mathcal{D}\colon W_1\to W_2$. We will write
$\langle \mathcal{D},w_1\rangle$ for the value of $\mathcal{D}$ on
$w_1\in W_1$. Let $M$ be an $H_D$-module. For all $m\in M$ we get an
$M$-valued distribution on $H_D$
\[
Y(m)\colon h\mapsto h.m,\quad h\in H_D.
\]
We can represent any distribution on $H_D$ by a generating series (or \emph{kernel}):
\[
\lambda(z)=\sum_{n=0}^\infty \langle\lambda,D^{(n)}\rangle z^n\in
M[[z]],\quad D^{(n)}=D^n/n!,
\]
and we have
\[
\lambda(h)=\epsilon(h.\lambda(z))=h.\lambda(z)|_{z=0},
\]
where $\epsilon\colon H_D^*\to \mathbb C$ is the counit and
where $h.$ indicates the action of $H_D$ on its dual $H_D^*$, with
$D.=\frac{d}{dz}$. In particular we have
\[
Y(m)(z)=e^{zD}m.
\]
In the theory of vertex algebras one usually writes $Y(m)(z)$ as
$Y(m,z)$. 

It is easy to write down a system of axioms for the \emph{holomorphic
  vertex operators} $Y(m,z)$ that is equivalent to the statement that
$M$ is a commutative, associative, unital algebra with a derivation.
In fact, one can easily define holomorphic vertex algebras for any
commutative and cocommutative Hopf-algebra $H$. These correspond to
commutative (etc.) rings with compatible $H$ action. In case $H$ is not
longer cocommutative a holomorphic vertex algebra corresponds to a
more general ring (braided for instance) with $H$-action. In this
paper we restrict ourselves to the commutative and cocommutative
Hopf algebras $H_D$ and $H_T$.

We call $Y(m,z)$ holomorphic since it has no singularities,  it
is a power series. The idea is now to define more general vertex
operators, that might contain singularities, while preserving the
axioms as much as is possible.

Introduce a singularization of the dual $H_D^*$:
\[
K_D=\mathbb C((z))=K_D^{\mathrm{Hol}}\oplus K_D^{\mathrm{Sing}},\quad
K_D^{\mathrm{Hol}}=\mathbb C[[z]],\quad K_D^{\mathrm{Sing}}=\mathbb
C[z\inv]z\inv.
\]
Then we have an isomorphism of $H_D$-modules:
\begin{equation}
  \label{eq:defalpha}
\alpha\colon H_D\to K_D^{\mathrm{Sing}}, \quad h\mapsto S(h)\frac1z.
\end{equation}
Here $S\colon H_D\to H_D$ is the antipode, where $D\mapsto -D$. Let
$\epsilon\colon H_D\to \mathbb C$ be the counit of $H_D$, the
multiplicative map such that $D\mapsto 0$. Observe then that
$\epsilon$ corresponds under $\alpha$ to the \emph{residue} on
$K_D^{\mathrm{Sing}}$.

We defined holomorphic vertex operators as distributions on
$H_D$. Using $\alpha$ we can think of them also as distributions on
$K_D^{\mathrm{Sing}}$. Then it is obvious how one can introduce
singular vertex operators: they should be distributions on all of
$K_D$, not just on $K_D^{\mathrm{Sing}}$.

Another ingredient in the theory of vertex algebras, besides vertex
operators and the residue, is the Dirac $\delta$-distribution, which
appears in the commutator of vertex operators. It can be described in
terms of the singularization $K_D$ of $H_D^*$ as follows. First of all
we have on $H_D^*$ the twisted coproduct $f(z)\mapsto f(z_1-z_2)$,
dual to the $S$-twisted multiplication $h_1\otimes h_2\mapsto
h_1S(h_2)$ on $H_D$. We can calculate this twisted coproduct in
\emph{two} ways, using $S$-twisted exponentials:
\begin{equation}
  \label{eq:twistedcoproduct}
f(z_1-z_2)=\mathcal{L}_S f(z_1)=\mathcal{R}_Sf(z_2),\quad
\mathcal{L}_S=e^{-z_2\partial_1},\quad \mathcal{R}_S=S_2
e^{z_1\partial_2}.
\end{equation}
If we try to extend the definition of the twisted coproduct from $H_D^*$ to
its singularization $K_D$ we see that the actions of $\mathcal{L}_S$
and $\mathcal{R}_S$ no longer give the same answer. We then define the
Dirac distribution for $p\in K_D$ as the obstruction to the extension
of the twisted coproduct:
\begin{equation}
\delta(p(z))=\mathcal{L}_S(
p(z_1))-\mathcal{R}_S(p(z_2)).\label{eq:defdeltap}
\end{equation}
For instance,
\[
\delta(\frac1z)=\delta(z_1,z_2)=\sum_{k\in\mathbb Z}z_1^kz_2^{-k-1}
\]
is the usual Dirac $\delta$-distribution.

Now we have formulated some ingredients of $H_D$-vertex algebras in
terms of the symmetry Hopf-algebra $H_D$ and the localization of the
dual, we are ready to try to extend the theory to other symmetry
algebras.

\section{$H_T$-symmetry and singularities.}
\label{sec:H_Tsym}

Let $H_T=\mathbb C[T,T\inv]$, and give it a Hopf algebra structure by
thinking of it as the group algebra of the free Abelian group generated
by $T$. This group is of course just the (additive) group of integers
$\mathbb Z$, and the dual of $H_T$ is the space
\[
H_T^*=\mathbb C_{\mathbb Z}=\{s\colon\mathbb Z\to \mathbb C\}
\]
of arbitrary complex valued functions on $\mathbb Z$, or,
equivalently, the space of two-sided infinite sequences
$s=(s_n)_{n\in\mathbb Z}$. We can expand any element $s\in\mathbb
C_{\mathbb Z}$ as an infinite sum
\[
s=\sum_{n\in\mathbb Z}s_n\delta_n,
\]
where $\delta_n$ is the Kronecker sequence; as a function on $\mathbb
Z$ we have $\delta_n\colon k\mapsto \delta_{kn}$. The natural action
of $H_T$ on its dual is given by $T\delta_n=\delta_{n-1}$.

Now we want to find a singularization of $H_T^*$. Note that it would
not be useful to invert some of the Kronecker sequences, as they are
zero divisors: $\delta_n\delta_k=\delta_{nk}\delta_n$. To find
suitable elements of $H_T^*$ to invert we observe that in the $H_D$ case
the powers $z^k$ we invert are the unique solutions of a system of
differential equations: if $f_\ell(z)\in H_D^*$ is a solution of
\[
\partial_z f_\ell(z)=\ell f_{\ell-1}(z),\quad f_\ell(z)|_{z=0}=0,
\quad \ell\ge 0,
\]
with $f_0(z)=1$, then $f_\ell(z)=z^\ell$.  Similarly, introduce in
$H_T$ the difference operator $\Delta=T-1$ and consider the system of
difference equations for $\tau(\ell)\in H_T^*$, $\ell\ge1$
\[
\Delta \tau(\ell)=\ell\tau(\ell-1), \quad\tau(\ell)|_0=0,
\]
where $\tau(0)=1=\sum_{n\in\mathbb Z}\delta_n$. Then the $\tau(\ell)$s
are uniquely determined. We have $\tau=\tau(1)=\sum_{n\in\mathbb
  Z}n\delta_n$ and $\tau(\ell)=\tau(\tau-1)\dots(\tau-\ell+1)$. The
$\tau(\ell)$s are the restriction of polynomial functions on $\mathbb
C$ to the integers.

Now let $\Czpol=\mathbb C[\tau]\subset H_T^*$, the space of polynomial
functions on $\mathbb Z$. Let then $M\subset \Czpol$ be the
multiplicative set generated by the translates $T^k\tau=\tau+k$ of
$\tau$. Then we define the following singularization 
\[
K_T=M\inv\Czpol.
\]
We could have localized all of $H_T^*$, but then there would have
been no clear distinction between singular and nonsingular elements in
the localization, see \cite{math.QA/0505289} for details.

Let $S:H_T\to H_T$ be the antipode, $T^k\mapsto T^{-k}$. We have, as in
the case of $H_D$, see (\ref{eq:defalpha}), a map
\[
\alpha\colon H_T\to K_T^{\mathrm{Sing}},\quad h\mapsto S(h)\frac1\tau,
\]
but it is not longer an isomorphism, for instance $\frac1{\tau^2}$ is
not in the image of $\alpha$. However, when we complete $H_T$
by adjoining to $H_T$ the infinite sum  
\[
\partial_\tau=\log(T)=\log(1+\Delta)=\sum_{n=1}^\infty(-\Delta)n/n,
\]
to define $\hat H_T=\mathbb C[T,T\inv,\partial_\tau]$, then
\[
\alpha\colon \hat H_T\to K_T^{\mathrm{Sing}},
\]
defined as before, \emph{is} an isomorphism. The counit $\epsilon\colon
H_T\to \mathbb C$ corresponds, via $\alpha$, to the map, called the \emph{trace},
\[
\Tr\colon K_T\to \mathbb C,\quad f(\tau)\mapsto \sum_{n\in\mathbb
  Z}\Res_{n}(f(\tau)d\tau).
\]
Associated to $H_T$ we have twisted exponential operators
$\mathcal{L}_S,\mathcal{R}_S$, analogous to those of
(\ref{eq:twistedcoproduct}), and we can define Dirac distributions as
before. In particular we have the Dirac $\delta$-distribution defined
by
\[
\delta(\frac1\tau)=\mathcal{L}_S\left(\frac1{\tau_1}\right)-
\mathcal{R}_S\left(\frac1{\tau_2}\right)=\sum_{n\in\mathbb
  Z}\tau(n)\otimes\tau(-n-1),\quad
\tau(-\abs{k})=\frac1{\tau(\abs{k})}.
\]
We have an action of $H_T\otimes H_T$ on such two-variable
distributions, and we write $h_1=h\otimes 1$, $h_2=1\otimes h$.
Similarly we write $\tau_1(\ell)=\tau(\ell)\otimes 1$ and
$\tau_2(\ell)=1\otimes \tau(\ell)$ and we denote the distribution
$\delta(\frac1\tau)$ also by $\delta(\tau_1,\tau_2)$. It has the usual
properties:
\begin{itemize}
\item $h_1\delta(\tau_1,\tau_2)=S(h)_2\delta(\tau_1,\tau_2)$, $h\in \hat
  H_T$.
\item $f(\tau_1)\delta(\tau_1,\tau_2)=f(\tau_2)\delta(\tau_1,\tau_2)$,
  $f(\tau)\in K_T$.
\item $\Tr_{\tau_1}(f(\tau_1)\delta(\tau_1,\tau_2))=f(\tau_2)$.
\item If a distribution $a(\tau_1,\tau_2)$ is killed by the twisted
  coproduct of an element $p\in\Czpol$ then it is a finite sum
  \begin{equation}
a(\tau_1,\tau_2)=\sum_{n,k}a_{n,k}(\tau_2)
e_{n,k}\delta(\tau_1,\tau_2),\label{eq:finiteexpansion}
\end{equation}
where $\{e_{n,k}\}$ is a basis for $\hat H_T$.
\end{itemize}

\section{Distributions and State-Field correspondence}
\label{sec:DistStateField}

Let $W$ be a vector space, and let $\mathcal{D}$ be a $W$-valued
distribution on $\KT$. Every such distribution $\mathcal{D}$ has a
formal expansion:
\[
\mathcal{D}(\tau)=\sum_{n\in\mathbb Z} \langle\mathcal{D},\tau(n)\rangle\tau(-n-1).
\]
If $F\in\KT$ then the value of $\mathcal{D}$ on $F$ is a trace:
\[
\langle\mathcal{D},F\rangle=\Tr(\mathcal{D}(\tau) F).
\]
We will often identify a distribution $\mathcal{D}$ with its kernel $\mathcal{D}(\tau)$.
Any distribution has a decomposition in holomorphic and singular part:
\[
\mathcal{D}=\mathcal{D}_{\Hol}+\mathcal{D}_{\Sing},
\]
where the kernels of the holomorphic and singular parts have expansion
in $\tau(n)$ for $n$ nonnegative, respectively negative. A
distribution $\mathcal{D}$ is called \emph{rational} if there is $\phi\in
W\otimes \KT$ such that $\langle \mathcal{D},F\rangle=\Tr(\phi F)$ for
all $F\in \KT$.

Now let $V$ be a vector space. A \emph{field} on $V$ is then an
$\End(V)$-valued distribution $f=f(\tau)$ such that for all $v\in V$
the $V$-valued distribution $f(\tau)v$ has rational singular part:
$f_{\Sing}(\tau)v\in V\otimes \KT^{\Sing}$. Denote by $V(\tau)$ the
space of fields on $V$. Then a \emph{State-Field Correspondence} is a
linear map $Y\colon V\to V(\tau)$. We write, if $a\in V$ and $Y$ is a
state-field correspondence, $Y(a)(\tau)=Y(a,\tau)$ and we call
$Y(a,\tau)$ the vertex operator of $a$, as usual. 

Now let $1=1_V\in V$ be a distinguished vector, called the
\emph{vacuum}. We say that a state-field correspondence $Y$
\emph{satisfies the vacuum axioms} (for $1\in V$) in case
\[
Y(1,\tau)=1_{\End(V)},\quad Y(f,\tau)1=f_{\Hol}(\tau)1,
\]
where $f_{\Hol}$ is a holomorphic distribution such that acting on the
vacuum the constant term is $f$: $f_{\Hol}(\tau)1|_{\tau=0}=f$.

Let $h\in H_T$. Then we define, given a state-field correspondence $Y$
satisfying the vacuum axioms, a linear map $h_V\colon V\to V$ by
\begin{equation}
h_Vf=\langle h, f_{\Hol}(\tau)\rangle 1,\label{eq:defh_V}
\end{equation}
where $\langle\,,\,\rangle$ is the pairing between $H_T$ and $H_T^*$,
extended in the obvious way to an $\End(V)$-valued pairing between
$H_T$ and holomorphic distributions. At this point we don't know that
the map $h\mapsto h_V$ gives an $H_T$-module structure to $V$, this is
an extra condition on the state-field correspondence.

If $a(\tau)$, $b(\tau)$ are fields on $V$, they are in particular
$\End(V)$-valued distributions, and we can calculate their
commutator distribution: this is the distribution (on $K_T\otimes
K_T$) that acts on $v\in
V$ by
\[
[a(\tau_1),b(\tau_2)](F\otimes G)v=\Big(a(F)b(G)-b(G)a(F)\Big)
v,\quad F,G\in K_T.
\]
We say that these fields are \emph{mutually rational} if the commutator
distribution $[a(\tau_{1}),b(\tau_{2})]$ has rational singularities,
i.e., is killed by some element $m_S^*(F)$, for $F\in H_T^*$ and
$m_S^*$ the twisted coproduct on $H_T^*$. In this case the commutator
is a finite sum of differential-differences of the delta distribution,
see (\ref{eq:finiteexpansion}).

\section{Definition and some properties of $H_T$-vertex algebras}
\label{sec:DefHTvertexAlg}

An $H_T$-vertex algebra is an $H_T$-module $V$ with a vacuum vector
$1_V\in V$ and a state-field correspondence $f\mapsto
f(\tau)=Y(f,\tau)$, satisfying the vacuum axioms and furthermore
\begin{itemize}
\item (Compatibility) The action of $H_T$ on $V$ is compatible with
  the state-field correspondence:
\[
h.f=h_V f,
\]
where the left hand side is the action of $H_T$ on $V$ and the right
hand side is defined in (\ref{eq:defh_V}).
\item ($\ad$-covariance) For all $f\in V$ and $h\in H_T$
\[
\ad^V_h Y(f, \tau)= h_{K_T}Y(f,\tau).
\]
Here $\ad^V_h(X)=\sum h^\prime X S(h^{\prime\prime})$, for $X\in\End(V)$,
$h\in H_T$ with coproduct $\Delta(h)=\sum h^\prime\otimes h^{\prime\prime}$.
\item (Mutual Rationality) The vertex operators $Y(f,\tau_1)$ and
  $Y(g,\tau_2)$ are for all $f,g\in V$ mutually rational.
\end{itemize}

From these axioms one derives easily properties similar to those of
$H_D$-vertex algebras.  For instance we have covariance of the
state-field correspondence:
\[
Y(h.f, \tau)=h_{K_T}Y(f,\tau),
\]
and skew-symmetry of vertex operators:
\[
Y(f,\tau)g=\mathcal{R}_V(\tau)Y^S(g,\tau)f,
\]
where $\mathcal{R}_V(\tau)$ is the exponential operator corresponding
to the $H_T$-action on $V$ and $Y^S$ is the antipodal vertex operator:
more generally if $\mathcal{D}$ is a distribution on $K_T$ then we
define its antipode by $\langle \mathcal{D}^S, F\rangle=\langle
\mathcal{D},S(F)\rangle$. 

We can define for $F\in K_T$ the $F$-product of $f,g$ in an
$H_T$-vertex algebra:
\[
f_{\{F\}}g=\Tr\Big(Y(f,\tau)gF(\tau)\Big)\in V.
\]
Also we can define the $F$-product of fields by
\[
f(\tau_2)_{\{F\}}g(\tau_2)=\Tr_{\tau_1}\Big(f(\tau_1)g(\tau_2)\mathcal{R}_S(\tau_2)F(\tau_1)-
g(\tau_2)f(\tau_1)\mathcal{L}_S(\tau_1)F(\tau_2)\Big).
\]
Here $\mathcal{R}_S,\mathcal{L}_S$ are the twisted exponentials
already used in (\ref{eq:defdeltap}). In particular, for
$F=\frac1\tau$ we obtain the \emph{normal ordered product} of fields:
\[
f(\tau)_{\{\frac1\tau\}}g(\tau)=
\nop{f(\tau)g(\tau)}=f_{\Hol}(\tau)g(\tau)-g(\tau)f_{\Sing}(\tau).
\]
Then we have the fundamental
fact that the state-field correspondence is a homomorphism of
$F$-products:
\[
Y(f_{\{F\}}g,\tau)=f(\tau)_{\{F\}}g(\tau).
\]

\section{$H_T$-conformal algebras}
\label{sec:H_Tconf}

For an $H_D$-vertex algebra $V_D$ we can concentrate on the singular
part of the operator expansion to obtain on $V_D$ the structure of
conformal algebra (see \cite{MR1651389}, or \cite{MR1670692}, where
conformal algebras are called vertex Lie algebras). In the same way we
can ignore in an $H_T$-vertex algebra $V$ all $F$-products
$f_{\{F\}}g$, except for those with $F\in\Czpol$, or more generally
$F\in H_T^*$. This defines the notion of an $H_T$-conformal algebra
structure on $V$.

More generally, one can start with an $H_T$-module $\mathcal{C}$ and
define an $H_T$-conformal structure on $\mathcal{C}$ as a collection
of \emph{conformal products} $f_{\{F\}}g$, $F\in H_T^*$ for $f,g\in
\mathcal{C}$, satisfying a number of axioms that we don't want write
down here\footnote{In fact, one can define for any cocommutative Hopf
  algebra an $H$-conformal algebra structure along these lines, see
  the notion of Lie pseudo algebra in Bakalov, D'Andrea and Kac,
  \cite{MR1849687}.}. In particular there is in an $H_T$-conformal
algebra the distinguished product corresponding to the element $F=1\in
H_T^*$. This product satisfies
\begin{align*}
  (Tf)_{\{1\}}g&= f_{\{1\}}g, & f_{\{1\}}(Tg)&= T(f_{\{1\}}g),\\
f_{\{1\}}g-g_{\{1\}}f&\in \mathfrak m_T\mathcal{C}, & [f_{\{1\}},g_{\{1\}}]&= (f_{\{1\}}g)_{\{1\}}.
\end{align*}
Here $\mathfrak m_T\subset H_T$ is the \emph{augmentation ideal}, the
kernel of the counit on $H_T$. It is the ideal generated by
$\Delta=T-1$. We see that the $1$-product induces a Lie algebra
structure on $\mathcal{C}/\mathfrak m_T \mathcal{C}$, for any
$H_T$-conformal algebra $\mathcal{C}$.

\section{Affinization}
\label{sec:affinization}
If $\mathcal{C}$ is an $H_T$-conformal algebra and $L$ is a
commutative $H_T$-Leibniz algebra (i.e., a commutative algebra in the
category of $H_T$-modules), then one can show that also the
\emph{affinization} $L\mathcal{C}=\mathcal{C}\otimes L$ is canonically
an $H_T$-conformal algebra, and hence we obtain on
$\mathcal{L}\mathcal{C}=L\mathcal{C}/\mathfrak m_T L\mathcal{C}$ the
structure of Lie algebra. In particular we can take $L=K_T$, and we
will restrict ourselves to this case. Denote by
\[
\Tr\colon \mathcal{C}\otimes K_T\to \mathcal{L}\mathcal{C}
\]
the canonical projection and write, for $p\in K_T$ and $f\in \mathcal{C}$,
\[
f_{\langle p \rangle}=\Tr(f\otimes p)\in \mathcal{L}\mathcal{C}.
\]
For an explanation of using the same term for both this map and the
trace on $\KT$, see \cite{math.QA/0505289}. Then the commutator in
$\mathcal{L}\mathcal{C}$ is given by
\begin{equation}
[ f_{\langle p \rangle},g_{\langle q \rangle}]=\sum (f_{\{
  e_i^*\}}g)_{\langle e_i(p)q\rangle},\label{eq:conformalcom}
\end{equation}
where $\{e_i\}$ is a basis for $H_T$ and $\{e_i^*\}$ a dual basis for
$H_T^*$. We can define generating series of elements of
$\mathcal{L}\mathcal{C}$, called \emph{currents}, for each $f\in \mathcal{C}$ by
\[
f(\tau_2)=\Tr_{\tau_1}\big(f\otimes \delta(\tau_1,\tau_2)\big).
\]
Then the commutator of currents is given by 
\begin{equation}
[f(\tau_1),g(\tau_2)]=\sum f_{\{e_i^*\}}g(\tau_2)
e_{i,2}\delta(\tau_1,\tau_2).\label{eq:currentcomm}
\end{equation}

\section{The Toda Vertex algebra}
\label{sec:Toda}

If $L$ is a commutative $H_T$-Leibniz algebra, then $L$ is
automatically an $H_T$-vertex algebra, giving the simplest examples
of them. In this case the vertex operators are holomorphic, see the
discussion in Section \ref{sec:RevH_DVertexalg}, so this is
not really interesting. 

To get a more interesting example we start with the Toda conformal
algebra $\CToda$. This is the free $H_T$-module generated by $B$ and
$C$, with conformal products 
\begin{align*}
  B_{\{F\}}B&=C_{\{F\}}C=0, \quad F\in H_T^* &\\
  B_{\{\delta_n\}}C&=C(\delta_{n,-1}-\delta_{n,0}), &
  C_{\{\delta_n\}}B&=C\delta_{n,0}-TC \delta_{n,1}.
\end{align*}
Here the $\delta_n\in H_T^*$ are the Kronecker sequences, see Section
\ref{sec:H_Tsym}. $\CToda$ is the $H_T$-conformal algebra
corresponding to the first Hamiltonian structure of the infinite Toda
lattice, see \cite{MR785802}. We write $\LToda$ for the Lie algebra
$\mathcal{L}\CToda$ associated to the Toda conformal algebra. 
The current commutator in $\LToda$ is
\begin{equation}
[B(\tau_1),C(\tau_2)]= C(\tau_2)
(T_2\inv-1)\delta(\tau_1,\tau_2).\label{eq:BCconformalcom}
\end{equation}
We have a decomposition
\begin{equation}
  \label{eq:LieTodaDecomp}
  \LToda=\LToda_{\Hol}\oplus \LToda_{\Sing},
\end{equation}
where $\LToda_{\Hol}$ (respectively $\LToda_{\Sing}$) is spanned by elements 
$f_{\langle p\rangle}$, for $p\in K_T^{\Hol}$ (respectively $p\in
K_T^{\Sing}$). By (\ref{eq:conformalcom}) each summand in
(\ref{eq:LieTodaDecomp}) is a Lie subalgebra. Let $\mathbb C_0$ be the
trivial 1-dimensional $\LToda_{\Hol}$-module, and let $\VToda$ be the
induced $\LToda$-module:
\[
\VToda=\mathcal{U}(\LToda)\otimes_{\mathcal{U}(\LToda_{\Hol})}\mathbb
C_0.
\]
Then one proves that $\VToda$ has an $H_T$-vertex algebra structure
(with $1_{\VToda}=1\otimes 1$ as vacuum) such that if we write $B$ and
$C$for the elements $B_{\langle \frac1\tau} 1_{\VToda}$ and
$C_{\langle \frac1\tau} 1_{\VToda}$, then we have
\[
Y(B,\tau_1),Y(C,\tau_2)]= Y(C,\tau_2)(T_2-1)\delta(\tau_1,\tau_2),
\]
compare with (\ref{eq:BCconformalcom}). We call this the Toda
$H_T$-vertex algebra structure on $\VToda$. 

\begin{remark}
    Note that, since $K_T$ is in fact an $\hat H_T$-module, also
  $\VToda$ will be not just an $H_T$-module, but an $\hat H_T$-module.
  In general, all non holomorphic $H_T$-vertex algebras seem to be
  $\hat H_T$-modules. This is in contrast with $H_T$-conformal
  algebras, which don't need to have an $\hat H_T$-modules structure.
  Now $\hat H_T=H_T[\partial_\tau]$ contains $H_D$, if we take
  $D=\partial_\tau$, and so we can think of nontrivial $H_T$-vertex
  algebras as a kind of extension of $H_D$-vertex algebras, where we
  allow in the operator product expansion not only singularities in
  $\tau_1-\tau_2$ but also at arbitrary shifts
  $T_2^k(\tau_1-\tau_2)=\tau_1-\tau_2-k$. Currents with such operator
  products expansions occur in the theory of Yangians, see e.g.,
  \cite{MR1601127}.
\end{remark}

\def\cprime{$'$}
\providecommand{\bysame}{\leavevmode\hbox to3em{\hrulefill}\thinspace}
\providecommand{\MR}{\relax\ifhmode\unskip\space\fi MR }
\providecommand{\MRhref}[2]{%
  \href{http://www.ams.org/mathscinet-getitem?mr=#1}{#2}
}
\providecommand{\href}[2]{#2}


\begin{thebibliography}{BDK01}

\bibitem[BDK01]{MR1849687}
Bojko Bakalov, Alessandro D'Andrea, and Victor~G. Kac, \emph{Theory of finite
  pseudoalgebras}, Adv. Math. \textbf{162} (2001), no.~1, 1--140. \MR{MR1849687
  (2003c:17020)}

\bibitem[Ber]{math.QA/0505289}
Maarten Bergvelt, \emph{{$H_T$-Vertex Algebras and the Infinite Toda Lattice}},
  arXiv:math.QA/0505289.

\bibitem[Bor86]{MR843307}
Richard~E. Borcherds, \emph{Vertex algebras, {K}ac-{M}oody algebras, and the
  {M}onster}, Proc. Nat. Acad. Sci. U.S.A. \textbf{83} (1986), no.~10,
  3068--3071. \MR{MR843307 (87m:17033)}

\bibitem[Bor98]{MR1653021}
\bysame, \emph{Vertex algebras}, Topological field theory, primitive forms and
  related topics (Kyoto, 1996), Progr. Math., vol. 160, Birkh\"auser Boston,
  Boston, MA, 1998, pp.~35--77. \MR{MR1653021 (99m:17034)}

\bibitem[Bor01]{MR1865087}
\bysame, \emph{Quantum vertex algebras}, Taniguchi Conference on Mathematics
  Nara '98, Adv. Stud. Pure Math., vol.~31, Math. Soc. Japan, Tokyo, 2001,
  pp.~51--74. \MR{MR1865087 (2002k:17054)}

\bibitem[FBZ04]{MR2082709}
Edward Frenkel and David Ben-Zvi, \emph{Vertex algebras and algebraic curves},
  second ed., Mathematical Surveys and Monographs, vol.~88, American
  Mathematical Society, Providence, RI, 2004. \MR{MR2082709 (2005d:17035)}

\bibitem[FLM88]{MR996026}
Igor Frenkel, James Lepowsky, and Arne Meurman, \emph{Vertex operator algebras
  and the {M}onster}, Pure and Applied Mathematics, vol. 134, Academic Press
  Inc., Boston, MA, 1988. \MR{MR996026 (90h:17026)}

\bibitem[Kac98]{MR1651389}
Victor Kac, \emph{Vertex algebras for beginners}, second ed., University
  Lecture Series, vol.~10, American Mathematical Society, Providence, RI, 1998.
  \MR{MR1651389 (99f:17033)}

\bibitem[Kho97]{MR1601127}
Sergej~M. Khoroshkin, \emph{Central extension of the {Y}angian double},
  Alg\`ebre non commutative, groupes quantiques et invariants (Reims, 1995),
  S\'emin. Congr., vol.~2, Soc. Math. France, Paris, 1997, pp.~119--135.
  \MR{MR1601127 (98j:17015)}

\bibitem[Kup85]{MR785802}
B.~A. Kuperschmidt, \emph{Discrete {L}ax equations and differential-difference
  calculus}, Ast\'erisque (1985), no.~123, 212. \MR{MR785802 (86m:58070)}

\bibitem[Pri99]{MR1670692}
Mirko Primc, \emph{Vertex algebras generated by {L}ie algebras}, J. Pure Appl.
  Algebra \textbf{135} (1999), no.~3, 253--293. \MR{MR1670692 (2000c:17046)}

\bibitem[Sny]{math.QA/9904104}
Craig~T. Snydal, \emph{{Equivalence of Borcherds G-Vertex Algebras and
  Axiomatic Vertex Algebras}}, arXiv:math.QA/9904104.

\end{thebibliography}
 \bibliographystyle{amsalpha}

\end{document}